\documentclass[a4paper]{amsart}
\usepackage[utf8]{inputenc}
\usepackage{multirow}
\usepackage{amssymb, amsmath}
\usepackage{amsfonts}
\usepackage{float}
\usepackage{tikz}
\usepackage{pgfplots}
\pgfplotsset{compat=1.12}
\usepackage{graphics}
\usepackage{subfigure}
\usepackage{xcolor}
\usepackage{graphicx}
\usepackage{enumitem}
\graphicspath{{./figures/}}

\newcommand{\abs}[1]{\lvert #1 \rvert}
\newcommand{\R}{\text{R}}
\newcommand{\divv}{\text{div}}
\newcommand{\Pre}{\text{Pr}}
\newcommand{\Le}{\text{Le}}

\newcommand{\np}[1]{d_{#1}}
\newcommand{\ip}[1]{\langle #1 \rangle}
\newcommand{\sign}{\operatorname{sign}}

\renewcommand{\Re}{\operatorname{Re}}

\newtheorem{lemma}{Lemma}
\newtheorem{theorem}{Theorem}

\begin{document}
\title[Transitions in Spherical THC]{Transitions of Spherical Thermohaline Circulation to Multiple Equilibria}

\author[Özer]{Saadet Özer}
\address[SO]{Department of Mathematics, Istanbul Technical University, 34469, Istanbul, Turkey}
\email{saadet.ozer@itu.edu.tr}

\author[Şengül]{Taylan Şengül}
\address[TS]{Department of Mathematics, Marmara University, 34722 Istanbul, Turkey}
\email{taylan.sengul@marmara.edu.tr}

\begin{abstract}
The main aim of the paper is to investigate the transitions of the thermohaline circulation in a spherical shell in a parameter regime which only allows transitions to multiple equilibria. We find that the first transition is either continuous (Type-I) or drastic (Type-II) depending on the sign of the transition number. The transition number depends on the system parameters and $l_c$, which is the common degree of spherical harmonics of the first critical eigenmodes, and it can be written as a sum of terms describing the nonlinear interactions of various modes with the critical modes. We obtain the exact formulas of this transition number for $l_c=1$ and $l_c=2$ cases. Numerically, we find that the main contribution to the transition number is due to nonlinear interactions with modes having zero wave number and the contribution from the nonlinear interactions with higher frequency modes is negligible. In our numerical experiments we encountered both types of transition for $\Le<1$ but only continuous transition for $\Le>1$. In the continuous transition scenario, we rigorously prove that an attractor in the phase space bifurcates which is homeomorphic to the 2$l_c$ dimensional sphere and consists entirely of degenerate steady state solutions.
\end{abstract}
\maketitle

\section{Introduction}
An important source of climate low frequency variability is the so called thermohaline circulation (THC). The underlying mechanism of THC is well established. THC is essentially driven by the temperature and freshwater fluxes in the ocean-atmosphere interface which in turn produces density gradients. These gradients are much sharper in the vertical direction compared to the horizontal directions and therefore are associated with an overturning.

There are indications that the Atlantic circulation has varied in the past \cite{broecker1991great}. Taking into account the enormous effects the ocean circulation has on the climate, the sensitivity, stability and transitions of the large scale ocean circulation became an important issue in climate research \cite{musundijkstra2004}. The studies of the THC using a hierarchy of ocean models, starting with a very simple box model \cite{stommel1961} shows that the presence of heat and salt, with their different influences on the density field, may lead to different stable steady flow patterns.

This paper arises out of a research program to generate rigorous mathematical results on climate variability developed from the viewpoint of dynamical transitions \cite{ptd,dijkstra2015dynamic,hsia2015tropical}. The basic philosophy of dynamic transition theory is to search for the full set of transition states, giving a complete characterization of stability and transition. The set of transition states which is often represented by a local attractor may lie near or away from the basic state. This theory has recently been successfully applied to many branches of nonlinear sciences, see \cite{yuan2016bifurcation, ong2016dynamic, choi2015bifurcation, yari2015transition} among others.

One of the main focus of the dynamic transition theory is the identification of transition states and the classification of dissipative systems into three distinct transition types, namely continuous (Type-I), catastrophic (Type-II) and random (Type-III) transitions which describe the nature of the state of the system as the control parameter crosses a critical threshold. The transition states stay in a close neighborhood of the basic state in the case of continuous transition, and lie outside of a neighborhood of the basic state in the case of a catastrophic transition. For a random type transition, a neighborhood of the basic state consists of two disjoint open regions with a continuous transition taking place in one region while a catastrophic transition occurring in the other one.

We briefly recall several recent papers which investigate the problem in simpler settings than the one studied in this paper. The paper \cite{bona2011hopf} considers the THC in a 2D rectangular box and considers the oscillatory transitions and finds that both continuous and drastic transitions are possible. In \cite{ma2010dynamic}, the authors investigate the dynamic transitions of THC in a 3D rectangular enclosure where they show that THC exhibits transitions to either multiple equilibria or to time periodic solutions and that there are parameter regimes leading to continuous or catastrophic transitions. Finally, in \cite{wang2013remarks}, the authors consider the dynamic transitions in a spherical shell but only for the pure thermal convection case without salinity gradients. The transition type in that case is well-known to be continuous and the transition is described by an attractor bifurcation \cite{ma2004dynamic}. Hence in \cite{wang2013remarks}, the authors aim to describe the structure of the local bifurcated attractor.

The main objective of this paper is to carry out the dynamic transition analysis of the THC problem in the important case of a spherical shell domain. For simplicity, the spatial domain is considered as the product of a two dimensional sphere $S_a^2$ with radius $a$ and an interval $(0 , h)$ where $h$ is is the height of the fluid layer. This approach is mainly motivated by the fact that the aspect ratio for the large scale atmosphere and ocean is small; see among others \cite{lions1992new,pedlosky2013geophysical}. The main challenge in carrying the analysis of \cite{ma2010dynamic} to a spherical shell domain is due to the detailed calculations of nonlinear interactions of spherical harmonics (both scalar and vectorial) that need to performed to compute the reduction of the infinite dimensional system to a finite one.

It is known \cite{ma2010dynamic} that the THC system exhibits first transitions to both multiple equilibria and to spatio-temporal oscillations depending on the parameter
\begin{equation} \label{K}
  K=\sign(1-\Le)\left[\frac{\Le^2}{1-\Le}(1+\frac{1}{\Pre})\sigma_c-\tilde{\R}\right],
\end{equation}
where $\Pre$ represents the Prandtl number, $\Le$ is the Lewis number, $\tilde{\R}$ is the saline Rayleigh number. Here $\sigma_c$ and $l_c \in \mathbb{Z}^+$ are defined by
\begin{equation} \label{sigma_c}
  \sigma_c = \min_{l \in \mathbb{Z}^+} \frac{(\pi^2 + \alpha_l^2)^3}{\alpha_l^2} = \frac{(\pi^2 + \alpha_{l_c}^2)^3}{\alpha_{l_c}^2},
\end{equation}
which exist thanks to the convexity of the function $x \to \frac{(\pi^2 + x)^3}{x}, \, x>0$. In \eqref{sigma_c}, the wave number $\alpha_l$ is defined as
\[
  \alpha_l^2 = \frac{l (l+1)}{r^2}, \quad l \in \mathbb{Z}^+,
\]
with $r=a/h$ denoting the aspect ratio of the sphere.
The well known result of Rayleigh Bénard convection is that the minimum of $\sigma_c$ is achieved when $\alpha_{l_c} = \frac{\pi}{\sqrt{2}}$ in which case $\sigma_c = \frac{27}{4} \pi^4$, see \cite{chandrasekhar1961}.

In this paper, we restrict ourselves to the parameter regime $K>0$ which only allows transitions to multiple equilibria and will address the case of transitions to spatio-temporal oscillations ($K<0$ regime) in another paper. In the case $K>0$, the first dynamic transition of the system occurs as the control parameter $\sigma$ defined by
\begin{equation} \label{sigma}
  \sigma = \R - \Le^{-1} \tilde{\R},
\end{equation}
crosses the critical threshold $\sigma_c$, leading to multiple equilibria. Here $\R$ is the thermal Rayleigh number.

Our main results are as follows. We find that the transition of the THC problem at $\sigma=\sigma_c$ for $l_c=1$ or $l_c=2$ is either Type-I or Type-II depending on the sign of the transition number $q_{l_c}$. We present the  exact formulas for the transition number $q_{l_c}$ which depend on the system parameters for $l_c=1$ and $l_c=2$ cases in \eqref{coefficientsq1}--\eqref{coefficientsq2}. In these cases, the transition number $q_{l_c}$ is given by the sum
\[
  q_{1} = \sum_{l=0}^{l_c} \np{(l_c,1),(2l,2)},
\]
where $\np{(l_c,1),(2l,2)}$ is a number determined by the nonlinear interactions of the critical modes with all modes corresponding to the spherical harmonics of degree $2l$, $0\le l \le l_c$ and vertical wave number $2$ (there are $2l+1$ of such modes).
In the continuous transition scenario, an attractor in the phase space bifurcates on $\sigma > \sigma_c$ which is a $2 l_c$-dimensional homological sphere. We show that this sphere is indeed homeomorphic to $S^{2l_c}$, the 2$l_c$ dimensional sphere and consists entirely of degenerate steady state solutions.

Our numerical explorations in the parameter space suggest that
\[
  \sign(q_{l_c}) =\sign(\np{(l_c,1),(0,2)}) = \frac{(1+\Le) \alpha_{l_c}^4 }{16 \pi \Le(\R_0 - \R)} (\R_1 - \R),
\]
when $\R$ is away from $\R_1 = \frac{\sigma_c}{1-\Le^2}$, indicating that the nonlinear interactions of the critical modes with the zero wave number modes determine the type of transition and that the contribution from higher frequency modes to the transition number is negligible. Also our numerical simulations suggest the following. For $\Le<1$, both continuous and drastic transitions are possible and a continuous transition is preferred for $\R < \R_{\ast}$ while a drastic transition is preferred for $\R_{\ast} < \R < \R_0 = \frac{\Le + \Pre}{(1-\Le) \Pre} \sigma_c$ with $R_{\ast} \approx R_1$. For $\Le>1$ flows, the transition type does not change and is always Type-I for all $\R$.

% For the large-scale atmospheric and oceanic circulations, the critical wave number $l_c$ is often relatively small. By introducing proper turbulent friction terms as in \cite{ma2010dynamic}, we can adjust the problem to match these convection scales. We shall report these issues elsewhere.

The paper is organized as follows. In Section 2, the mathematical setting of the problem is introduced. In Section 3, the linear stability analysis of the main equations is summarized. Section 4 details our main theorem and its proof. Finally in Section 5, we explore the transition number numerically.

\section{Governing Equations and the Functional Setting}
As mentioned in the Introduction, we consider a spherical shell $S_a^2\times (0,h)$ as the spatial domain for the motion of the large scale ocean. Here $S_a^2$ represents the 2D sphere with radius $a$ and $h$ denotes the height of the fluid layer. The governing equations are the familiar Boussinesq equations (see \cite{pedlosky2013geophysical,ma2010dynamic} among others):
\begin{equation} \label{field1}
\begin{aligned}
  &\frac{\partial {\bf u}}{\partial t}+({\bf u} \cdot \nabla){\bf u}=\nu\Delta {\bf u}-\frac{1}{\rho_0} (\nabla p+\rho g  \hat{e}_z),\\
  &\frac{\partial T}{\partial t}+({\bf u} \cdot \nabla )T=\kappa_T\Delta T,\\
  &\frac{\partial S}{\partial t}+({\bf u} \cdot \nabla )S=\kappa_S\Delta S,\\
  &\divv {\bf u}=0,
\end{aligned}
\end{equation}
where ${\bf u}$ is the velocity, $T$ is the temperature, $S$ is the salinity, $\hat{e}_z$ is the unit vector in the z-direction, $\nu$, $\kappa_T$, $\kappa_S$, $g$ are all positive constants denoting the kinematic diffusivity, the thermal diffusivity, the saline diffusivity and the gravitational constant, respectively. The fluid density $\rho$ is given by the linear equation of state
\begin{equation} \label{rho}
  \rho=\rho_0 \left (1-a_T(T-T_0)+a_S(S-S_0) \right),
\end{equation}
where $a_T$ and $a_S$ are assumed to be positive constants and $\rho_0$ is the density at the lower surface. The constants $T_0$, $S_0$ represent the fixed temperature and salinity at the lower boundary $z=0$, whereas we denote the fixed temperature and salinity at the upper boundary $z=h$ by $T_1$ and $S_1$.

In \eqref{rho}, the opposing effects of temperature and salinity are evident. Both $T_0>T_1$ (heated from below) and $S_0<S_1$ (salted from above) are destabilizing mechanisms. We will treat the general case which allows one or both of these conditions to be satisfied. In other words, a competition between a stabilizing and a destabilizing mechanism is allowed.

The trivial steady state solution to the problem \eqref{field1}--\eqref{rho} is given by
\begin{equation}\label{steady state solutions}
\begin{aligned}
  &{\bf u}_{\text{ss}}=0,\\
  &T_{\text{ss}}=T_0-(T_0-T_1)\frac{z}{h},\\
  &S_{\text{ss}}=S_0-(S_0-S_1)\frac{z}{h},\\
  &p_{\text{ss}}=p_0-g\rho_0 \left( z+\frac{a_T}{2}(T_0-T_1)\frac{z^2}{h}-\frac{a_S}{2}(S_0-S_1)\frac{z^2}{h} \right).
\end{aligned}
\end{equation}
We nondimensionalize the equations \eqref{field1} exactly in the same way as in \cite{ma2010dynamic}, and obtain

% To nondimensionalize the equation we consider the perturbation of the solution from the steady state
% \begin{align*}
% &u''=u-u^0,&T''=T-T^0\\
% &S''=S-S^0,&p''=p-p^0
% \end{align*} and we set
% \begin{align*}
% &x=hx',&t=h^2t'/ \kappa_T\\
% &u''=\kappa_T u'/h,&T''=(T_0-T_1)/T'\\
% &S''=\abs{S_0-S_1}S',&p''=\rho_0\nu \kappa_T p'/h^2
% \end{align*}
% Omitting the primes the equation, \eqref{field1} can be written as

\begin{equation} \label{field2}
\begin{aligned}
  &\frac{\partial {\bf u}}{\partial t}=\Pre(\Delta {\bf u}-\nabla p)+\Pre(\R T-\sign(S_0-S_1)\tilde{\R}S) \hat{e}_z-({\bf u} \cdot \nabla){\bf u},\\
  &\frac{\partial T}{\partial t}=\Delta T+w-({\bf u} \cdot \nabla )T,\\
  &\frac{\partial S}{\partial t}=\Le \Delta S+\sign(S_0-S_1)w-({\bf u} \cdot \nabla )S,\\
  &\divv {\bf u}=0,
\end{aligned}
\end{equation}
where $w$ is the vertical velocity in the direction of $\hat{e}_z$, $\R$ is the thermal Rayleigh number, $\tilde{\R}$ is the saline Rayleigh number, $\Pre$ is the Prandtl number and $\Le$ is the Lewis number defined as
\[
  \begin{aligned}
    & \R=\frac{a_T g(T_0-T_1)h^3}{\kappa_T\nu}, &&\tilde{\R}=\frac{a_S g(S_0-S_1)h^3}{\kappa_T\nu}, \\
    & \Pre=\frac{\nu}{\kappa_T}, &&\Le=\frac{\kappa_S}{\kappa_T}.
  \end{aligned}
\]

In \eqref{field2}, the unknowns now represent deviations from the steady state solutions given by \eqref{steady state solutions}.
Also the nondimensional spatial domain is $\Omega=S_r^2\times (0,1)$ where the aspect ratio is
\begin{equation} \label{aspect ratio r}
  r=\frac{a}{h}.
\end{equation}
Let ${\bf u} = (u,w)$ be the 3D velocity vector where $u = u_{\theta} \hat{e}_{\theta} + u_{\varphi} \hat{e}_{\varphi}$ is the 2D horizontal velocity field. Hereafter, $\nabla$, $\nabla_u $, div and $\Delta$ will denote both the scalar and vectorial differential operators in the horizontal direction, that is on the sphere $S^2_r$, given by
\[
\begin{aligned}
  &\nabla_u f = u \cdot \nabla f = \frac{1}{r}\left(u_{\theta}\frac{\partial f}{\partial\theta}+\frac{u_{\varphi}}{\sin\theta}\frac{\partial f}{\partial\varphi}\right), \\
  &\Delta f=\frac{1}{r^2\sin\theta}\left[\frac{\partial}{\partial\theta}\left(\sin\theta\frac{\partial f}{\partial\theta}\right)+\frac{1}{\sin\theta}\frac{\partial^2f}{\partial\varphi^2}\right], \\
  &\nabla_u v=\frac{1}{r}\left(u_{\theta}\frac{\partial v_{\theta}}{\partial \theta}+\frac{u_{\varphi}}{\sin \theta}\frac{\partial v_{\theta}}{\partial \varphi}-u_{\varphi}v_{\varphi}\cot\theta\right)\hat{e}_{\theta}
  +\frac{1}{r}\left(u_{\theta}\frac{\partial v_{\varphi}}
  {\partial \theta}+\frac{u_{\varphi}}{\sin\theta}
  \frac{\partial v_{\varphi}}{\partial\varphi}+u_{\varphi}
  v_{\theta}\cot\theta\right)\hat{e}_{\varphi}, \\
  &\Delta u=\left(\Delta u_{\theta}-\frac{2\cos\theta}{r^2\sin^2\theta}\frac{\partial u_{\varphi}}{\partial \varphi}-\frac{u_{\theta}}{r^2\sin^2\theta}\right)\hat{e}_{\theta}+\left(\Delta u_{\varphi}+\frac{2\cos\theta}{r^2\sin^2\theta}\frac{\partial u_{\theta}}{\partial \varphi}-\frac{u_{\varphi}}{r^2\sin^2\theta}\right)\hat{e}_{\varphi},
\end{aligned}
\]
where $f$ is a scalar, $u=u_{\theta} \hat{e}_{\theta} + u_{\varphi} \hat{e}_{\varphi}$ and $v=v_{\theta} \hat{e}_{\theta} + v_{\varphi} \hat{e}_{\varphi}$ are 2D vectors.

With the above notations, \eqref{field2} can be written as
\begin{equation} \label{field3}
\begin{aligned}
  &u_t+\nabla_u u+w  \frac{\partial u}{\partial z}+\Pre\nabla p-\Pre(\Delta+\partial_{zz})u=0, \\
  &w_t+\nabla_u w+w \frac{\partial w}{\partial z}+\Pre \  \frac{\partial p}{\partial z}-\Pre(\R T-\sign(S_0-S_1)\tilde{\R}S)-\Pre(\Delta+\partial_{zz})w=0,\\
  &T_t+\nabla_u T+w \frac{\partial T}{\partial z}-w-(\Delta+\partial_{zz})T=0,\\
  &S_t+\nabla_u S+w \frac{\partial S}{\partial z}-\sign(S_0-S_1)w-\Le(\Delta+\partial_{zz})S=0, \\
  &\divv u+ \frac{\partial w}{\partial z}=0.
\end{aligned}
\end{equation}
In this study, we consider the equations \eqref{field3} supplemented with the free-slip boundary conditions
\begin{equation} \label{bc}
w=T=S=0, \quad \frac{\partial u}{\partial z}=0, \quad \text{at } z=0,1,
\end{equation}
and note that our analysis can be expanded to other types of boundary conditions as well.

The below functional spaces are needed to recast the equations \eqref{field3} and \eqref{bc} in an abstract form.
\begin{align*}
&H=\{(u,w,T,S)\in L^2(\Omega)^5\mid \divv u = w=T=S=0 {\text{ at }}  z=0,1\},\\
&H_1=\{(u,w,T,S)\in H^2(\Omega)^5\mid \divv u= w=T=S=0, \, \frac{\partial u}{\partial z} =0, \text{ at } \ z=0,1\}.
\end{align*}
We recall that the inner product in $H$ is defined for vectors $\Psi_i = (u_i, w_i, T_i, S_i)$ as
\[
\ip{\Psi_1, \Psi_2} = \int_0^{\pi}\int_0^{2\pi}\int_0^1 (u_1 \cdot \overline{u_2} + w_1 \overline{w_2} + T_1 \overline{T_2} + S_1 \overline{S_2}) r^2 \sin \theta dz d\varphi d\theta.
\]
Let the linear operator $L :H_1\rightarrow H$ be defined by
\[
  L (\Psi)= \mathcal{P}
\begin{bmatrix}
  \Pre(\Delta+\frac{\partial^2}{\partial z^2})u \\
  \Pre(\Delta+\frac{\partial^2}{\partial z^2})w+\Pre(\R T-\sign(S_0-S_1)\tilde{\R}S) \\
  (\Delta+\frac{\partial^2}{\partial z^2})T+w \\
  \Le(\Delta+\frac{\partial^2}{\partial z^2})S+\sign(S_0-S_1)w
\end{bmatrix},
\quad \Psi = (u, w, T, S) \in H_1,
\]
and the bilinear operator $G:H_1 \times H_1 \rightarrow H$ be defined by
\[
G(\Psi_1, \Psi_2) = -\mathcal{P}
\begin{bmatrix}
  \nabla_{u_1} u_2 + w_1 \frac{\partial u_2}{\partial z} \\
  \nabla_{u_1} w_2 + w_1 \frac{\partial w_2}{\partial z} \\
  \nabla_{u_1} T_2 + w_1 \frac{\partial T_2}{\partial z} \\
  \nabla_{u_1} S_2 + w_1 \frac{\partial S_2}{\partial z}
\end{bmatrix}, \quad
\Psi_i = (u_i, w_i, T_i, S_i) \in H_1.
\]
Here $\mathcal{P}:L^2(\Omega)^5\rightarrow H$ is the Leray projection.

Now the problem can be cast as an abstract ODE as
\begin{equation} \label{abstract equ}
\frac{\partial \Psi}{\partial t}=L \Psi+G(\Psi),\qquad \Psi(0)=\Psi_0,
\end{equation}
where $G(\Psi)=G(\Psi,\Psi).$

\section{Linear Stability}
We first consider the eigenvalue problem for the linearized equations of \eqref{field3}
\begin{align}
  &\Pre \left( (\Delta + \partial_{zz}) u - \nabla p \right)=\beta u \label{linearu1},\\
  &\Pre \left( (\Delta + \partial_{zz}) w + \R T-\sign(S_0-S_1)\tilde{\R}S - \frac{\partial p}{\partial z} \right) = \beta w, \label{linearw1} \\
  &(\Delta+\partial_{zz})T+w=\beta T, \nonumber \\
  &\Le(\Delta+\partial_{zz})S+\sign(S_0-S_1)w=\beta S, \nonumber \\
  &\divv u+\frac{\partial w}{\partial z}=0, \label{linearinc}
\end{align}
supplemented with the boundary conditions with \eqref{bc}.

The z-independent solutions of \eqref{linearu1}--\eqref{linearinc} clearly satisfy $w=T=S=0$.
In this case, taking the curl of \eqref{linearu1} and using the incompressibility condition is equivalent to
\[
  \Pr \nabla \times (\nabla \times (\nabla \times u)) = \beta ( \nabla \times u).
\]
It is easy to obtain the solutions of this equation and they are given by \eqref{beta_l0}.

% Using the fact that, see \cite{barrera85},
% \[
%   \nabla \times \nabla \times \nabla \times \nabla \times {\bf Y}_{lm} = \alpha_l^2 \nabla \times \nabla \times {\bf Y}_{lm}
% \]
% where ${\bf{ Y}}_{lm} = Y_{lm} \hat{e}_z$, we find a family of eigenpairs given by \eqref{beta_l0}.
% \[
%   u = \nabla \times {\bf Y}_{lm} = \frac{1}{r} \frac{1}{\sin \theta} \frac{\partial Y_{lm}}{\partial \varphi} \hat{e}_{\theta} - \frac{1}{r} \frac{\partial Y_{lm}}{\partial \theta} \hat{e}_{\varphi}
% \]
% \url{http://notes.yeshiwei.com/_static/Vector_spherical_harmonics.pdf}

For the general case, we use the separation of variables in the form
\begin{equation}\label{separation}
\begin{aligned}
  & u=\nabla f(\theta,\varphi)H'(z), \\
  & w=\alpha^2 f(\theta,\varphi)H(z), \\
  & T=f(\theta,\varphi)\Theta(z), \\
  & S=f(\theta,\varphi)\Phi(z).
\end{aligned}
\end{equation}
With \eqref{separation}, the incompressibility equation  \eqref{linearinc} is equivalent to the Helmholtz equation on the sphere
\begin{equation} \label{helmholtz}
  \Delta f+\alpha^2 f=0,
\end{equation}
which has solutions $f=Y_{lm}(\theta,\varphi)$ only when $\alpha^2 = \alpha_l^2 = \frac{l(l+1)}{r^2}$, where $Y_{lm}$ are the spherical harmonics, $l,m \in \mathbb{Z}$, $l \ge 1$ and $-l \le m \le l$.

Eliminating the pressure term $p$ by taking $\nabla$\eqref{linearw1}-$\frac{\partial}{\partial z}$\eqref{linearu1},  using \eqref{separation}--\eqref{helmholtz} and denoting $D = \frac{d}{dz}$, the equations \eqref{linearu1}--\eqref{linearinc} become
\begin{equation} \label{linearfi2}
\Pre \left( (D^2 - \alpha^2)^2 H -  \R \Theta + \sign(S_0-S_1)\tilde{\R}\Phi \right) = \beta (D^2 - \alpha^2)H,
\end{equation}
\begin{equation} \label{lineart2}
(D^2 -\alpha^2) \Theta + \alpha^2 H = \beta\Theta,
\end{equation}
\begin{equation} \label{linears2}
\Le(D^2 - \alpha^2) \Phi + \sign(S_0-S_1)\alpha^2 H = \beta \Phi.
\end{equation}
From \eqref{bc}, the boundary conditions of the above system of ODEs are
\[
  (\Theta,\Phi,H, D^2 H)=0,\ {\text {at}} \ z=0,1.
\]

If $l = 0$, then $u = 0$, $w=0$ and the equations \eqref{linearfi2}--\eqref{linears2} reduce to
\begin{equation*}
\begin{aligned}
  & D^2 \Theta = \beta \Theta, \\
  & \Le D^2 \Phi = \beta\Phi,
\end{aligned}
\end{equation*}
whose solutions are given in \eqref{beta_0n}.

On the other hand, if $l \ne 0$, plugging
\begin{align*}
& \Theta = \Theta_{lmn} \sin n\pi z, \\
& \Phi = \Phi_{lmn} \sin n\pi z, \\
& H = \sin n \pi z,
\end{align*}
into \eqref{linearfi2}, \eqref{lineart2} and \eqref{linears2}, we find the coefficients $\Theta_{lmn}$ and $\Phi_{lmn}$ given by \eqref{Psi lmn} and the compatibility condition \eqref{dispersion equ} for the existence of eigenvalues.

\subsection{Eigenpairs}
We will denote the eigenvalues of the linearized operator by $\beta_{ln}$ and their corresponding eigenvectors by $\Psi_{lmn}=(u, w, T, S)$ where $l,m,n \in \mathbb{Z}$, $l\ge 0$, $-l\le m \le l$, $n \ge 0$. By the discussion in the previous section, we now summarize all the eigenpairs.

If $l \ne 0$ and $n=0$, then
\begin{equation} \label{beta_l0}
  \beta_{l0}=-\Pre\alpha_l^2, \quad \Psi_{lm0}=(\nabla \times Y_{lm} \hat{e}_z,0,0,0), \qquad \abs{m} \le l,
\end{equation}
where
\[
  \nabla \times Y_{lm}\hat{e}_z = \frac{1}{r \sin \theta}\frac{\partial Y_{lm}}{\partial \varphi} \hat{e}_{\theta}-\frac{1}{r}\frac{\partial Y_{lm}}{\partial\theta}\hat{e}_{\varphi}.
\]

If $l = 0$ and $n\ne 0$, then
\begin{equation} \label{beta_0n}
  \begin{aligned}
    & \beta_{0n}^1=-n^2\pi^2, && \Psi_{00n}^1=(u=0,w=0,T=\sin n\pi z,S=0),\\
    & \beta_{0n}^2=-\Le\ n^2\pi^2, && \Psi_{00n}^2=(u=0,w=0,T=0,S=\sin n\pi z).
  \end{aligned}
\end{equation}

If $l \ne 0$ and $n \ne 0$, there are three distinct eigenvalues which we order as $\Re(\beta^1_{ln})\geq \Re(\beta^2_{ln})\geq \Re(\beta^3_{ln})$ corresponding to the three distinct solutions of
\begin{equation} \label{dispersion equ}
\beta^3 + b_2 \beta^2 + b_1 \beta + b_0 =0,
\end{equation}
where
\[
  \begin{aligned}
    & b_0 = (n^2\pi^2+\alpha_l^2)^3 \Le \Pre - \alpha_l^2\Pre( \Le  \R - \tilde{\R}), \\
    & b_1 = (n^2\pi^2+\alpha_l^2)^2 (\Le + \Pre + \Le \Pre) - \alpha_l^2(n^2\pi^2+\alpha_l^2)^{-1} \Pre( \R -\tilde{\R}), \\
    & b_2 = (n^2\pi^2+\alpha_l^2)(1 + \Le +  \Pre).
  \end{aligned}
\]

As the relation \eqref{dispersion equ} is independent of $m$, to each $\beta_{ln}^k$ ($l \ne 0$, $n \ne 0$, $k\in\{1,2,3\}$) corresponds $2 l + 1$ eigenfunctions $\Psi_{lmn}^k$, $m = -l, \dots, l$ given by
\begin{equation} \label{Psi lmn}
\Psi_{lmn}^k=
\begin{cases}
  & u= n\pi \nabla Y_{lm} \cos n\pi z, \\
  & w=\alpha_l^2 Y_{lm} \sin n\pi z, \\
  & T=\dfrac{\alpha_l^2}{n^2\pi^2+\alpha_l^2+\beta_{ln}^k}Y_{lm} \sin n\pi z, \\
  & S=\dfrac{\alpha_l^2 \sign(S_0-S_1)}{\Le(n^2\pi^2+\alpha_l^2)+\beta_{ln}^k}Y_{lm} \sin n\pi z.
\end{cases}
\end{equation}

By the symmetry $\overline{Y}_{lm} = (-1)^m Y_{l-m}$ of the spherical harmonics, it follows that if $\overline{\beta}_{ln}^{k_1} = \beta_{ln}^{k_2}$ then
\begin{equation} \label{symmetry of eigens}
  \overline{\Psi}_{lmn}^{k_1} = (-1)^m \Psi_{l-mn}^{k_2}.
\end{equation}
In particular when $\beta_{ln}^k \in \mathbb{R}$, we have $\overline{\Psi}_{lmn}^{k} = (-1)^m \Psi_{l-mn}^{k}$.

\subsection{Adjoint Problem}
For our analysis, we also require the explicit form of the adjoint eigenvectors. The adjoint problem of \eqref{linearu1}-\eqref{linearinc} can be written as
\begin{align*}
&\Pre\left((\Delta+\partial_{zz})u^*-\nabla p^*\right)=\overline{\beta} u^* , \\
&\Pre \left((\Delta+\partial_{zz})w^*-\frac{\partial p}{\partial z}^*\right)+ T^*-\sign(S_0-S_1)S^*=\overline{\beta}w^*,\\
&(\Delta+\partial_{zz})T^*+\Pre\R w^*=\overline{\beta}T^*, \\
&\Le(\Delta+\partial_{zz})S^*-\Pre \tilde{\R}\sign(S_0-S_1)w^*=\overline{\beta}S^*, \\
&\divv u^*+\frac{\partial w}{\partial z}^*=0.
\end{align*}
In the case $n=0$, $l \ne 0$ and in the case $l=0$, $n \ne 0$, the adjoint eigenfunctions are as given in \eqref{beta_l0} and \eqref{beta_0n} respectively.
When $l\neq 0$ and $n \ne 0$, the eigenfunctions are now given by
\[
  \Psi_{lmn}^{k*}=
  \begin{cases}
    & u=n\pi \nabla Y_{lm} \cos n\pi z,\\
    & w=\alpha_l^2 Y_{lm} \sin n\pi z,\\
    & T=\dfrac{\Pre \R\alpha_l^2}{n^2\pi^2+\alpha_l^2+\overline{\beta}_{ln}^k}Y_{lm} \sin n\pi z,\\
    & S=-\dfrac{\Pre\tilde{\R}\alpha_l^2 \sign(S_0-S_1)}{\Le(n^2\pi^2+\alpha_l^2)+\overline{\beta}_{ln}^k}Y_{lm} \sin n\pi z.
  \end{cases}
\]
\subsection{Principle of Exchange of Stabilities}
A crucial parameter in the first transition is the critical wave integer $l_c$ defined in \eqref{sigma_c} which depends only on the aspect ratio $r$. This dependence can be easily obtained explicitly from \eqref{sigma_c} by solving $\sigma_c(l) = \sigma_c(l+1)$ and the following lemma holds.
\begin{lemma}
  Define
  \begin{equation} \label{rl}
    r_l^2 = \frac{(1+l)}{\pi^{2}}   \left( \sqrt[3]{l(2+l)^2} + \sqrt[3]{l^2 (2+l)} \right), \quad l = 0,1,2,\dots
  \end{equation}
  Then the critical integer $l_c$ defined implicitly by \eqref{sigma_c} is
  \begin{equation} \label{l_c}
    l_c = l, \qquad \text{if } r_{l-1} < r < r_l,
  \end{equation}
  % Asymptotic Behaviour: Moreover, $r_l = \frac{\sqrt{2}}{\pi} (l+1) +O(l^{-1})$ as $l \to \infty$ so that if $r \gg 1$ and $(l-1) \le \frac{\pi}{\sqrt{2}} r \le l$, then $l_c = l$.
\end{lemma}
Numerical values of the first three $r_l$, computed from \eqref{rl}, are
\[
  r_0 = 0, \qquad r_1 = 0.844851, \qquad r_2 = 1.31566.
\]

We remark here that when the aspect ratio is critical, i.e. $r=r_l$ for some $l \in \mathbb{Z}^+$, the eigenvalues $\beta_{l_c1}^1$ and $\beta_{l_c+1,1}^1$ with wave numbers $\alpha_{l_c}$ and $\alpha_{l_c+1}$ will become critical simultaneously. This type of transition was studied for the pure 2D Bénard convection in a rectangular container in \cite{sengul-rb2d} leading to competition among the pure $l_c$-modes, the pure $l_c+1$-modes and the mixed modes (superpositions of $l_c$ and $l_c+1$ modes).

The dispersion relation \eqref{dispersion equ} is exactly the same as in \cite{ma2010dynamic} in which the thermohaline circulation problem is considered in a rectangular domain.
Thus the following principle of exchange of stabilities follows directly from the corresponding one in \cite{ma2010dynamic}.
\begin{theorem} \label{TheoremPES}
Assume that $\Le \ne 1$ and that the aspect ratio $r$ defined by \eqref{aspect ratio r} is not critical, i.e. $r \ne r_l$ for any $l \in \mathbb{Z}^+$, with $r_l$ as defined by \eqref{rl}.
Consider $K$, $\sigma$, $\sigma_c$ and the integer $l_c$ defined by \eqref{K}, \eqref{sigma}, \eqref{sigma_c} and \eqref{l_c} respectively.

\begin{enumerate}
\item If $K>0$ then $\sigma_c$ is the first critical Rayleigh number and the first critical eigenvalue $\beta^1_{l_c1}$ has multiplicity $2l_c + 1$. Moreover the condition
\[
\begin{aligned}
  & \beta^1_{l_c1}
    \begin{cases}
      < 0, & \sigma < \sigma_c ,\\
      = 0, & \sigma = \sigma_c ,\\
      > 0, & \sigma > \sigma_c,
    \end{cases} \\
  & \Re(\beta^k_{ln}) < 0, \qquad \text{if } (k, l, n) \ne (1, l_c, 1),
\end{aligned}
\]
is satisfied.

\item If $K<0$, let
\[
  \eta = R - \frac{\Pre + \Le}{\Pre + 1} \tilde{R}.
\]
Then
\[
  \eta_c = \frac{(\Pre + \Le)(1+\Le)}{\Pre} \sigma_c,
\]
is the first critical Rayleigh number. Moreover the first two critical eigenvalues satisfy $\beta^1_{l_c1} = \overline{\beta^2_{l_c1}}$ and each has multiplicity $2l_c + 1$ satisfying
\[
  \begin{aligned}
    & \Re \beta^1_{l_c1} = \Re \beta^2_{l_c1}
      \begin{cases}
        < 0, & \eta < \eta_c ,\\
        = 0, & \eta = \eta_c ,\\
        > 0, & \eta > \eta_c,
      \end{cases} \\
    & \Re(\beta^k_{ln}) < 0, \qquad \text{if } (k, l, n) \ne (1, l_c, 1) \text{ or } (k, l, n) \ne (2, l_c, 1).
  \end{aligned}
\]
\end{enumerate}
\end{theorem}
% In this paper we limit ourselves to the study of $K>0$ case.
% In the case $K<0$, the first transition leads to spatio-temporal oscillations.
Since the sign of the parameter $K$ defined by \eqref{K} plays a crucial role in the transition, let us investigate it in some detail. For this let us define
\[
  \R_0 = \frac{\Le + \Pre}{(1-\Le) \Pre} \sigma_c.
\]
Then
\[
  K=
  \begin{cases}
    \Le (\R_0 - \R) \sign(1-\Le), & \text{when } \sigma = \sigma_c, \\
    \dfrac{\Pre+1}{\Pre + \Le} (\R_0 - \R) \sign(1-\Le), & \text{when } \eta = \eta_c.
  \end{cases}
\]

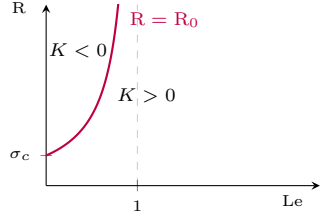
\begin{figure}
\centering
\pgfplotsset{every tick label/.append style={font=\tiny}, every axis/.append style={font=\tiny}}
\begin{tikzpicture}
  \begin{axis}[
  axis lines=middle, % left, right, box, center, none
  xmin=0, xmax=3, ymin=0, ymax=6,
  x=12mm, y=4mm, % width of the image
  xlabel={$\Le$},
  ylabel= {$\R$},
  x label style={at={(axis description cs:0.9,0)},anchor=north},
  y label style={at={(axis description cs:-0.1,.9)},anchor=south}, % rotate=90,
  xtick={1},
  ytick={1}, % ytick=\empty,
  xticklabels={1}, % xticklabels={,,},
  yticklabels={$\sigma_c$},
  ]
  \addplot[domain=0:.99, thick, purple, samples=100] {(x+3)/(3*(1-x)) };
  \addplot[domain=0:6, dashed, opacity=.2] ({1},{x});
  \node[] at (axis cs: .35, 4.5) {\scriptsize $K<0$};
  \node[] at (axis cs: 1.1, 3) {\scriptsize $K>0$};
  \node[purple] at (axis cs:1.3,5.5) [] {\scriptsize $\R = \R_0$};
\end{axis}
\end{tikzpicture}
\caption{The sign of the parameter $K$ which distinguishes the oscillatory and steady convection regimes in the $\Le-\R$ plane . \label{Fig: sign of K}}
\end{figure}

If $\Le>1$, it is clear that $\R_0<0$ and $K>0$.
From these observations, we find that
\begin{equation} \label{signOfK}
  K
    \begin{cases}
    > 0 , & \text{if } \Le >1 \text{ or } \R < \R_0, \\
    < 0 , & \text{if } \Le <1 \text{ and } \R > \R_0,
    \end{cases}
\end{equation}
at the criticality $\sigma=\sigma_c$ or $\eta=\eta_c$, see also Figure~\ref{Fig: sign of K}.

\section{Main Theorem and its Proof}

In this section we describe the transitions of the system \eqref{abstract equ} to multiple equilibria. For simplicity, we will consider the cases for which the critical integer $l_c$ defined by \eqref{sigma_c} is $1$ or $2$. The reason for this assumption is that the reduction to the center manifold, although still computable, becomes increasingly difficult to compute as $l_c$ is increased.

We recall that according to Lemma~\ref{sigma_c}, the aspect ratio ratio $r=a/h$ must be in the range $0<r<r_1$ for the $l_c=1$ case and $r_1<r<r_2$ for the $l_c=2$ case.

Now let us define the following coefficients
\begin{equation} \label{transition number coefficients}
  \begin{aligned}
    & a_{l}^k= \Le(\alpha_l^2 +4\pi^2) + \beta_{l2}^k,\\
    & b_{l}^k= \alpha_l^2+4\pi^2+ \beta_{l2}^k, \\
    & c_{l}^k = \frac{\alpha_{l_c}^4}{(\pi^2+\alpha_{l_c}^2)^2}\left( \sigma + (\pi^2 + \alpha_{l_c}^2) \Pre \left( \frac{\R}{b_{l}^k} - \frac{\tilde{\R}}{a_{l}^k \Le} \right) \right), \\
    & f_{l}^k = 4 \pi^2 + \alpha_l^2 \left( 1 + \Pre \left( \frac{\R}{(b_{l}^k)^2} - \frac{\tilde{\R}}{(a_{l}^k)^2}\right)\right), \\
    & g = \frac{\alpha_{l_c}^2}{(\pi^2+\alpha_{l_c}^2)^ 2} \frac{\Pre (1-\Le)}{\Le} (\R_0 - \R), \quad
    \R_0 = \frac{\Le + \Pre}{(1-\Le) \Pre} \sigma_c,
  \end{aligned}
\end{equation}

\begin{equation}\label{D02l}
  \np{(l_c,1),(0,2)} = \frac{(1+\Le) \alpha_{l_c}^4 }{16 \pi \Le(\R_0 - \R)} (\R_1 - \R), \quad
  \R_1 = \frac{\sigma_c}{1-\Le^2},
\end{equation}

\begin{equation} \label{coefficientsq1}
  \begin{aligned}
    & q_1 = \np{(1,1),(0,2)} + \np{(1,1),(2,2)}, \\
    & \np{(1,1),(2,2)} = \frac{3\pi}{40\sqrt{2}} \frac{1}{g}\sum_{k=1}^3 \frac{(c_{2}^k)^2}{\beta_{22}^k f_{2}^k},
  \end{aligned}
\end{equation}

\begin{equation} \label{coefficientsq2}
  \begin{aligned}
    & q_2 = \np{(2,1),(0,2)}+\np{(2,1),(2,2)}+\np{(2,1),(4,2)}, \\
    & \np{(2,1),(2,2)} = \frac{45 \pi }{784} \frac{1}{g}
    \sum_{k=1}^3 \frac{(c_{2}^k)^2}{\beta_{22}^k f_{2}^k}, \\
    & \np{(2,1),(4,2)} = -\frac{5 \pi}{42} \frac{1}{g}
    \sum_{k=1}^3 \frac{(c_{4}^k)^2}{\beta_{42}^k f_{4}^k}.
  \end{aligned}
\end{equation}
In the above notation, $\np{(l_c, 1),(l, n)}$ is the coefficient arising from the nonlinear interactions of the critical modes, i.e. eigenmodes corresponding to the $\beta^1_{l_c1}$ eigenvalue with all the modes corresponding to $\beta_{ln}$ eigenvalue.

% The term $\np{(l_c, 1),(0, 2)}$ measuring the nonlinear interaction of the critical modes with all modes corresponding to eigenvalue $\beta_{0,2}$ for all integers $l_c \ge 1$ is computed in \eqref{D02l}. The higher frequency contributions are only computed for $l_c=1$ and $l_c=2$ and are given in \eqref{coefficientsq1} and \eqref{coefficientsq2}.

\begin{theorem}
  Assume the conditions of Theorem~\ref{TheoremPES} and that $l_c=1$ or $l_c=2$. Consider $q_{l_c}$ defined by \eqref{coefficientsq1} and \eqref{coefficientsq2}. If $\sigma < \sigma_c$ then the basic solution $\psi = 0 \in H$ is locally asymptotically stable.
  \begin{itemize}
    \item If $q_{l_c} > 0$ then there is a continuous (Type-I) transition at $\sigma = \sigma_c$ and an attractor $\Sigma_{\sigma}$ bifurcates on $\sigma > \sigma_c$ which is homeomorphic to the 2$l_c$ dimensional sphere. $\Sigma_{\sigma}$ consists of degenerate steady states and has the following approximation
    \[
      \Sigma_{\sigma} = \left\{ \sum_{m=-l_c}^{l_c} x_{m} \Psi_{l_c m 1}^1 \mid \sum_{m=-l_c}^{l_c} \abs{x_m}^2 = \frac{\beta_{l_c1}^1 (\sigma)}{q_{l_c}} + O(\abs{\beta_{l_c1}^1}) \right \}.
    \]
    Moreover there is an open neighborhood $U \in H$ of $0$ such that $\Sigma_{\sigma}$ attracts $U \setminus \Gamma_{\sigma}$ where $\Gamma_{\sigma}$ is the stable set of $0$ of codimension $2l_c + 1$.
    \item If $q_{l_c} < 0$ then there is a drastic (Type-II) transition at $\sigma = \sigma_c$ and a repeller bifurcates on $\sigma < \sigma_c$ and no steady state solution bifurcates on $\sigma > \sigma_c$.
  \end{itemize}
\end{theorem}

\subsection{Proof of the Main Theorem}
The proof lies mainly on the reduction of the equation \eqref{abstract equ} onto the center manifold which is tangent to the critical eigenspace
\[
  E_1= \text{span} \left \{ \sum_{m=-l_c}^{l_c} x_m \Psi^1_{l_c m 1} \in H \mid x_{-m} = (-1)^m \overline{x_m}, \, x_m \in \mathbb{C} \right \},
\]
near the criticality $\sigma=\sigma_c$ and a careful analysis of the reduced equations. The condition $x_{-m} = (-1)^m \overline{x_m}$ is required in light of \eqref{symmetry of eigens} so that $E_1$ consists of real valued functions.

Let
\[
  \Phi_{c} = \sum_{m=-l_c}^{l_c} x_m(t) \Psi^1_{l_c m 1} \in E_1.
\]

Now a crucial element of the reduction procedure is the determination of the approximation of the center manifold function which captures the local dynamics near the criticality $\sigma = \sigma_c$. We will make use of the following approximation of the center manifold, see \cite{ptd},
\begin{equation} \label{cm approximation formula}
  -\mathcal{L} \Phi_{cm} = P_2 G(\Phi_{c})+o(2),
\end{equation}
in our analysis. Here $P_2:H\rightarrow E_2 := \{\Psi \in H_1 \mid \ip{\Psi, \Psi^{1 \ast}_{l_c,m,1}}=0 \text{ for all } \abs{m}\le l_c\}$ is the canonical projection onto the orthogonal complement, $\mathcal{L}$ is the restriction of the linear operator $L$ in \eqref{abstract equ} onto $E_2$ and
\[
  o(2) = o\left( \abs{x}^2 \right) + O(\abs{\beta^1_{l_c1}(\sigma)}\abs{x}^2), \quad \text{as } \sigma \to \sigma_c, \, \abs{x} \to 0,
\]
where $\abs{x}^2 = \sum_{m=-l_c}^{l_c} \abs{x_m}^2$.

According to \eqref{cm approximation formula}, the center manifold is determined by the following nonlinear interactions of the critical modes
\[
  \ip{ G(\Phi_{c}),\Psi^{\ast}_{lmn} } = \sum_{m_1=-l_c}^{l_c} \sum_{m_2=-l_c}^{l_c} x_{m_1}x_{m_2} \ip{ G \left( \Psi_{l_c m_1 1}^1, \Psi_{l_c m_2 1}^1 \right),\Psi^{\ast}_{lmn} },
\]
where
\begin{equation} \label{nonlin approx explicit}
  \begin{aligned}
    & \ip{ G \left( \Psi_1, \Psi_2 \right),\Psi^{\ast}_{lmn} } =
    \int_{\Omega} \bigg[(\nabla_{u_1} u_2 + w_1 \frac{\partial u_2}{\partial z}) \cdot(n\pi \nabla \overline{Y_{lm}} \cos n\pi z)\\
    & +( \nabla_{u_1} w_2 + w_1 \frac{\partial w_2}{\partial z})\alpha_l^2 \overline{Y_{lm}} \sin n\pi z  \\
    & +(\nabla_{u_1} T_2 + w_1 \frac{\partial T_2}{\partial z})(\frac{\Pre \R\alpha_l^2}{n^2\pi^2+\alpha_l^2+\beta_{ln}} \overline{Y_{lm}} \sin n\pi z)\\
    & +(\nabla_{u_1} S_2 + w_1 \frac{\partial S_2}{\partial z})(\frac{\Pre\tilde{\R}\alpha_l^2 \sign(S_0-S_1)}{\Le(n^2\pi^2+\alpha_l^2)+ \beta_{ln}} \overline{Y_{lm}} \sin n\pi z) \bigg] r^2 \sin \theta dz d\varphi d\theta,
  \end{aligned}
\end{equation}
for $\Psi_j = (u_j, w_j, T_j, S_j)$, $j=1,2$.

Due to \eqref{nonlin approx explicit}, the integral of triple product of spherical harmonics
\[
  \int_0^{\pi} \int_{0}^{2\pi} Y_{l_1m_1} Y_{l_2 m_2} \overline{Y}_{l_3 m_3} d\Omega = c_{l_1,l_2,l_3, m_1, m_2, m_3},
\]
plays a crucial role in the determination of the above nonlinear interactions. The coefficients $c_{l_1,l_2,l_3, m_1, m_2, m_3}$ are related to Clebsch-Gordan coefficients and are zero if
\[
  m_1 + m_2 \ne m_3,
\]
or
\[
  \abs{l_{j_1} - l_{j_2}} > l_{j_3} \text{ for some distinct } j_1, j_2, j_3 \in \{1,2,3\}.
\]
Since the critical modes always have wave index 1 in the z-direction, the product in \eqref{nonlin approx explicit} will vanish unless $n=0$ or $n=2$. As a summary of the above remarks, we obtain the following expansion of the center manifold function
\begin{equation} \label{Phi}
  \Phi_{cm} = \sum_{1 \le l \le 2 l_c,\abs{m}\leq l} y_{lm0}\Psi_{lm0}+ \sum_{k=1}^2 y_{002}^k \Psi_{002}^k + \sum_{\substack{1 \le l \le 2 l_c, \, \abs{m}\leq l, \\ 1\le k \le 3}} y_{lm2}^k\Psi_{lm2}^k + o(2),
\end{equation}
where by \eqref{cm approximation formula},
\begin{equation} \label{cm coefficients formula}
  \begin{aligned}
    & y_{lm0}=-\frac{\langle G(\Phi_{c}),\Psi_{lm0}^{\ast}\rangle}{\beta_{l0}\ip{\Psi_{lm0}, \Psi_{lm0}^{\ast}}}, \\
    & y_{002}^k=-\frac{ \ip{ G(\Phi_{c}),\Psi_{002}^{k\ast} } }{\beta_{02}^k\ip{\Psi_{002}^{k}, \Psi_{002}^{k\ast}}}, \\
    & y_{lmn}^k=-\frac{\langle G(\Phi_{c}),\Psi_{lmn}^{k\ast}\rangle}{\beta_{ln}^k\ip{\Psi_{lmn}^{k}, \Psi_{lmn}^{k\ast}}}.
  \end{aligned}
\end{equation}
Next we carry out the computation of the center manifold coefficients in \eqref{cm coefficients formula} for $l_c = 1$ and $l_c = 2$ separately. In both cases, we find
\begin{equation} \label{ylm0}
  y_{lm0} = 0,
\end{equation}
for all $1 \le l \le 2 l_c$ and $\abs{m} \le l_c$.

From \eqref{Phi} and \eqref{ylm0}, we may write the center manifold for $l_c=1$ as
\[
  \Phi_{cm} = \sum_{k=1}^2y_{002}^k\Psi_{002}^k + \sum_{k=1}^3\sum_{m=-2}^2 y_{2m2}^k\Psi_{2m2}^k,
\]
where,
\begin{align*}
& y_{002}^k=A_{02}^k(x_0^2-2x_{-1}x_1), \quad k=1,2, \\
& \frac{y_{2-22}^k}{x_{-1}^2} = \frac{y_{222}^k}{x_{1}^2} = \frac{y_{2-12}^k}{\sqrt{2}x_{-1}x_0} = \frac{y_{212}^k}{\sqrt{2}x_{1}x_0} = \frac{y_{202}^k}{\sqrt{\frac{2}{3}}(x_0^2+x_{-1}x_1)} = A_{22}^k ,\quad k=1,2,3, \\
& A_{02}^1 = \frac{ \alpha_{l_c}^4 }{ 16 \pi^2 ( \pi^2+\alpha_{lc}^2)},\quad
A_{02}^2=\frac{ \sign(S_0-S_1)}{\Le^2}A_{02}^1, \quad
A_{22}^k = -\frac{1}{4} \sqrt{\frac{3\pi}{10}} \frac{c_{2}^k}{\beta_{22}^k f_{2}^k},
\end{align*}
and $f_{l}^k$, $c_{l}^k$ are given by \eqref{transition number coefficients}.

From \eqref{Phi} and \eqref{ylm0}, we may write the center manifold for $l_c=2$ as
\[
  \Phi_{cm}=\sum_{k=1}^2y_{002}^k\Psi_{002}^k+\sum_{k=1}^3\sum_{m=-2}^2y_{2m2}^k\Psi_{2m2}^k+\sum_{k=1}^3\sum_{m=-4}^4y_{4m2}^k\Psi_{4m2}^k.
\]
Here
\begin{align*}
    y_{002}^k & =A_{02}^k(x_0^2-2x_{-1}x_1+2x_{-2}x_2), \quad k=1,2, \\
    \begin{split}
    B_{22}^k & = \frac{y_{2-22}^k}{(4 x_0x_{-2}-\sqrt{6}x_{-1}^2)} =
    \frac{y_{222}^k}{(4x_0x_2-\sqrt{6}x_1^2)} =
    \frac{y_{2-12}^k}{2(\sqrt{6}x_1x_{-2}-x_0x_{-1})} \\
    & = \frac{y_{212}^k}{2(\sqrt{6}x_2x_{-1}-x_0x_1)} = \frac{y_{202}^k}{2(2x_2x_{-2}+x_1x_{-1}-x_0^2)}
  \end{split}, \quad k=1,2,3,\\
  \begin{split}
    B_{42}^k & =
    \frac{y_{4-42}^k}{x_{-2}^2} =
    \frac{y_{4-32}^k}{\sqrt{2}x_{-2}x_{-1}} =
    \frac{y_{4-22}^k}{\frac{1}{\sqrt{7}}(2x_{-1}^2+\sqrt{6}x_0x_{-2})} \\
    & = \frac{y_{4-12}^k}{\sqrt{\frac{2}{7}} (x_1x_{-2}+\sqrt{6}x_0x_{-1})} =
    \frac{y_{402}^k}{\sqrt{\frac{2}{35}} (4x_1x_{-1}+x_2x_{-2}+3x_0^2)} \\
    & = \frac{y_{412}^k}{\sqrt{\frac{2}{7}} (x_1x_2+\sqrt{6}x_0x_1)} =
    \frac{y_{422}^k}{\frac{1}{\sqrt{7}}(2x_1^2+\sqrt{6}x_0x_2)} =
    \frac{y_{432}^k}{\sqrt{2}x_2x_1} =
    \frac{y_{442}^k}{x_{2}^2}
  \end{split}, \quad k=1,2,3, \\
   B_{22}^k & = \frac{3\sqrt{5\pi}}{56}\frac{c_{2}^k}{\beta_{22}^k f_{2}^k}, \quad
  B_{42}^k = -\frac{1}{6}\sqrt{\frac{5\pi}{14}} \frac{c_{4}^k}{\beta_{42}^k f_{4}^k},
\end{align*}
and $f_{l}^k$, $c_{l}^k$ are given by \eqref{transition number coefficients}.

Now, we are ready to write down the reduced equations. We plug in $\Psi = \Phi_{c} + \Phi_{cm}$ into \eqref{abstract equ} and project the resulting equation onto $E_1$ to obtain
\begin{equation}\label{reducedmain}
\frac{dx_m}{dt}=\beta_{l_c1}^1\ x_m+\frac{\langle G(\Phi_{c}+\Phi_{cm}),\Psi_{l_cm1}^{1\ast}\rangle}{\langle\Psi_{l_c m1}^1, \Psi_{l_c m1}^{1\ast}\rangle} + o(3), \qquad m = -l_c, \dots, l_c.
\end{equation}
Since $\int_0^1 s_1(z) s_2(z) s_3(z) dz = 0$ where $s_i(z) = \cos \pi z$ or $s_i(z) = \sin \pi z$, one gets right away that
\[
  \ip{G(\Phi_{c}, \Phi_{c}), \Psi_{l_c m 1}^{1 \ast}} = 0.
\]
Also as $\Phi_{cm} = O(\abs{x}^2)$, one gets $G(\Phi_{cm}, \Phi_{cm}) = O(\abs{x}^4)$ and hence the equations \eqref{reducedmain} become
\begin{equation}\label{reducedmain2}
\frac{dx_m}{dt}=\beta_{l_c1}^1\ x_m+\frac{ \ip{ G(\Phi_{c}, \Phi_{cm})  + G(\Phi_{cm}, \Phi_{c}),\Psi_{l_cm1}^{1\ast} }  }{\langle\Psi_{l_c m1}^1, \Psi_{l_c m1}^{1\ast}\rangle} + o(3),
\end{equation}
for each $m \in \{-l_c, \dots, l_c\}$.

We compactly present the results of tedious computations, both manual and by symbolic computation software, of the nonlinear terms in \eqref{reducedmain2} for the cases $l_c = 1$ and $l_c = 2$ separately below. For the $l_c = 1$ case, the reduced equations become
\begin{equation} \label{reduced equation lc=1}
  \frac{dx_m}{dt}= \beta^1_{11}\ x_m - q_1 x_m(x_0^2-2x_{-1}x_1)+o(3),\quad m=-1,0,1,
\end{equation}
where $q_1$ is defined as in \eqref{coefficientsq1}.

For the $l_c=2$ case, the reduced equations are
\begin{equation} \label{reduced equation lc=2}
\frac{dx_{m}}{dt}=\beta^1_{21}\ x_m - q_2 x_m(x_0^2-2x_{-1}x_1+2x_{-2}x_2) + o(3),\qquad m=-2,-1,0,1,2,
\end{equation}
where $q_2$ is defined as in \eqref{coefficientsq2}.

Since $x_{-m} x_m = (-1)^m \frac{1}{2} (\abs{x_{-m}}^2 + \abs{x_m}^2)$, letting
\[
  \abs{x}^2 = \sum_{m=-l_c}^{l_c} \abs{x_m}^2,
\]
the reduced equations \eqref{reduced equation lc=1} and \eqref{reduced equation lc=2} can be written compactly as
\begin{equation} \label{reduced equation final}
\frac{dx_m}{dt}=\beta^1_{l_c1}\ x_m - q_{l_c} x_m \abs{x}^2 + o(3),\quad m=-l_c,\dots, l_c.
\end{equation}
By taking the product of \eqref{reduced equation final} with $\overline{x_m}$ and summing from $m=-l_c$ to $m=l_c$, we derive
\begin{equation} \label{reduced equation amplitude form}
  \frac{1}{2}\frac{d \abs{x}^2}{dt} = \beta^1_{l_c1} \abs{x}^2 - q_{l_c} \abs{x}^4 + o(4).
\end{equation}
Since the critical eigenspace is $2l_c+1$ dimensional which is odd, using Krasnosel'skii Theorem (see Theorem 1.10 in \cite{ma-wang2005}) the existence of a bifurcated nontrivial steady state solution of the main equations at $\sigma =\sigma_c$ can be shown exactly as in \cite{wang2013remarks}.

In the case $q_{l_c}>0$ and $\sigma > \sigma_c$ case, from \eqref{reduced equation amplitude form}, it follows that $\abs{x}^2 \to \frac{\beta^{1}_{l_c1}(\sigma)}{q_{l_c}}$ as $t \to \infty$ for all sufficiently small initial conditions. By the attractor bifurcation theorem in \cite{ma-wang2005}, the bifurcated attractor $\Sigma_{\sigma}$ is homeomorphic to the $2l_c$ dimensional sphere $S^{2l_c}$. Moreover since the main equations possesses $S^{2l_c}$-symmetry, this steady state solution will generate a $S^{2l_c}$ set of steady states. Thus $\Sigma_{\sigma}$ consists solely of steady state solutions which are all degenerate since the Jacobian determinant of the right hand side of \eqref{reduced equation final} vanishes at the criticality $\beta^1_{l_c1}(\sigma)=0$ for those steady state solutions. That proves the assertions of our main theorem.

% The stability analysis of \eqref{reduced equation amplitude form} is easy for sufficiently small $\abs{x}$. If $q_{l_c} > 0$ then $\abs{x}^2 \to \frac{\beta_{l_c1}^1}{q_{l_c}}$ when $\beta^1_{l_c1} > 0$, i.e. $\sigma > \sigma_c$ as $t \to \infty$. Hence in this case, $2l_c$ dimensional sphere of steady states bifurcate on $\sigma > \sigma_c$.

\section{Numerical Computations of the Transition Number}
As shown in our main theorem, the type of first transition to multiple equilibria depends on a nondimensional parameter $q_{l_c}$, $l_c \in \mathbb{Z}^+$ which in turn depends on five system parameters: the aspect ratio $r$, the Lewis number $\Le$, the Prandtl number $\Pre$, the thermal Rayleigh number $\R$, the saline Rayleigh number $\tilde{\R}$. These parameters are related by the equation $\R- \Le^{-1} \tilde{\R} = \sigma_c$, where $\sigma_c$ is given in \eqref{sigma_c} at the onset of transition, leaving four degrees of freedom in the determination of $q_{l_c}$. Hence, in the following discussion we fix $\tilde{\R}$ by the choice of other four parameters. Also the parameter regime we are interested in our main theorem is the region $K>0$ where $K$ is given by \eqref{K}. According to \eqref{signOfK} this corresponds to either $\Le >1$ regime or $\Le <1$ and $\R < \R_0 = \dfrac{\Le + \Pre}{(1-\Le) \Pre} \sigma_c$ regime, see Figure~\ref{Fig: sign of K}.

Our main theorems for $l_c = 1$ and $l_c = 2$ cases prove that the type of transition is governed by the parameters $q_1 = \np{(1,1),(0,2)} + \np{(1,1),(2,2)}$ and $q_2 = \np{(2,1),(0,2)} + \np{(2,1),(2,2)} + \np{(2,1),(4,2)}$, respectively. We recall that $\np{(l_c, 1),(l, n)}$ is the coefficient arising from the nonlinear interactions of the critical modes with all the modes corresponding to the $\beta_{ln}$ eigenvalue.

We first present the numerical values of these nonlinear interaction terms when $\Pre=7.5$,  $\Le=10^{-2}$, $\Le=10^{-1}$, $\Le=0.5$, $\Le=5$ and $\R = 620$, $\R=640$ for the case $l_c=1$ (with $r=2/\pi$) in Table~\ref{tbl: D1102 D1122} and $l_c = 2$ (with $r=2\sqrt{3}/\pi$) in Table~\ref{tbl: D2102 D2122 D2142}.
\begin{table}[ht]
\centering
\caption{The numerical values of the nonlinear interaction terms for $l_c=1$, $r=2/\pi$, $\Pre=7.5 $ for different $\Le$ and $\R$.
\label{tbl: D1102 D1122}
}
\begin{tabular}{|l|l|l|l|l|}
\hline
& \multicolumn{2}{|c|}{$\R=620$} & \multicolumn{2}{|c|}{$\R=660$} \\
\hline
$\Le$     & $\np{(1,1),(0,2)}$ & $\np{(1,1),(2,2)}$ & $\np{(1,1),(0,2)}$ & $\np{(1,1),(2,2)}$ \\ \hline
$10^{-2}$ & 40.825     & 1.493      & -23.53   & -0.8537    \\ \hline
$10^{-1}$ & 1.955    & 0.074      & 0.275   & 0.0122     \\ \hline
0.5       & 0.477     & 0.0196     & 0.424   & 0.0176     \\ \hline
5         & 0.421	& 0.175			& 0.427	& 0.1177	\\   \hline
\end{tabular}
\end{table}

\begin{table}[ht]
\centering
\caption{The numerical values of the nonlinear interaction terms for $l_c=2$, $r=2 \sqrt{3}/\pi$, $\Pre=7.5 $ for different $\Le$ and $\R$.
\label{tbl: D2102 D2122 D2142}
}
\begin{tabular}{|l|l|l|l|l|l|l|}
\hline
& \multicolumn{3}{|c|}{$\R=620$} & \multicolumn{3}{|c|}{$\R=660$} \\
\hline
$\Le$     & $\np{(2,1),(0,2)}$ & $\np{(2,1),(2,2)}$ & $\np{(2,1),(4,2)}$ & $\np{(2,1),(0,2)}$ & $\np{(2,1),(2,2)}$ & $\np{(2,1),(4,2)}$ \\ \hline
$10^{-2}$ & 40.825     & 8.359   & 0.592    &-23.53 & -4.797 & -0.338   \\ \hline
$10^{-1}$ & 1.956    & 0.407   & 0.029   & 0.275  & 0.063  & 0.005   \\ \hline
0.5       & 0.477     & 0.104  & 0.008    & 0.424   & 0.093  & 0.007   \\ \hline
5		  & 0.421	& 0.092	& 0.0069		  & 0.427	& 0.093    & 0.007  \\ \hline
\end{tabular}
\end{table}

\begin{figure}
\centering
\pgfplotsset{every tick label/.append style={font=\tiny}, every axis/.append style={font=\tiny}}
\begin{tikzpicture}
  \begin{axis}[
  axis lines=middle, % left, right, box, center, none
  xmin=0, xmax=3, ymin=0, ymax=6,
  x=14mm, y=4mm, % width of the image
  xlabel={$\Le$},
  ylabel= {$\R$},
  x label style={at={(axis description cs:0.9,0)},anchor=north},
  y label style={at={(axis description cs:-0.1,.9)},anchor=south}, % rotate=90,
  xtick={2},
  ytick={1}, % ytick=\empty,
  xticklabels={1}, % xticklabels={,,},
  yticklabels={$\sigma_c$},
  ]
  \addplot[domain=0:.75, thick, purple, samples=100] {(x+3)/(3*(1-x)) };
  \addplot[domain=0:1.7, thick, blue, samples=100] {(x+4)/(4-x^2) };
  \addplot[domain=0:6, dashed, opacity=.2] ({2},{x});
  \node[purple] at (axis cs:1,5.5) [] {\scriptsize $\R = \R_0$};
  \node[blue] at (axis cs:2,5.5) [] {\scriptsize $\R = \R_1$};
  \node[] at (axis cs:2.3,2.5) [] {\scriptsize $+$};
  \node[] at (axis cs:1.6,2.0) [] {\scriptsize $+$};
  \node[] at (axis cs:1.1,3.5) [] {\scriptsize $-$};
\end{axis}
\end{tikzpicture}
\caption{The sign of $\np{(l_c,1),(0,2)}$ in the $\Le-\R$ plane. \label{Di02}}
\end{figure}
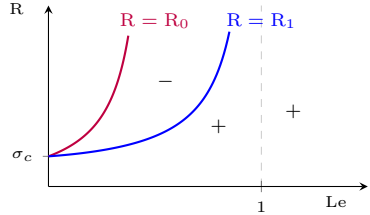

Our numerical results in Table~\ref{tbl: D1102 D1122} and Table~\ref{tbl: D2102 D2122 D2142} strongly indicate that
\[
  \frac{\abs{\np{(1,1),(0,2)}}}{\abs{\np{(1,1),(2,2)}}} \gg 1, \qquad
  \frac{\abs{\np{(2,1),(0,2)}}}{\abs{\np{(2,1),(2,2)}}+\abs{\np{(2,1),(4,2)}}} \gg 1.
\]
Thus the sign of $\np{(l_c,1),(0,2)}$ given by \eqref{D02l} plays a crucial role in the determination of the transition to multiple equilibria. So we investigate the sign of $\np{(l_c,1),(0,2)}$ in detail now. Recalling that $\R_1 = \frac{\sigma_c}{1-\Le^2}$, we observe that $\R_0 > \R_1 > 0$ if $\Le < 1$ and $\R_1<0$ if $\Le > 1$. Thus we obtain, see also Figure~\ref{Di02},
\[
  \sign(\np{(l_c,1),(0,2)}) =
  \begin{cases}
    -\sign(\R-\R_1), & \Le < 1, \\
    1, & \Le > 1.
  \end{cases}
\]

The above remarks suggest that $\sign(q_1) = \sign(\np{(1,1),(0,2)})$ and $\sign(q_2) = \sign(\np{(2,1),(0,2)})$ at least when $\R$ is away from $\R_1 = \frac{\sigma_c}{1-\Le^2}$.
Hence the contribution to the transition number $q_{l_c}$ of the nonlinear interactions of the critical modes with the higher frequency modes is negligible compared to the interactions with zero wavenumber modes.

\begin{figure}[ht]
\centering
% Le = 10^{-2},
\begin{tikzpicture}
\pgfplotsset{every tick label/.append style={font=\tiny}, every axis/.append style={font=\tiny}, }
\begin{axis}[
xlabel=$\R$,
ylabel = $q_i$,
% scaled ticks=false,
axis lines=center,
y label style={at={(axis description cs:-0.1,.9)},anchor=south},
title={$\Le=10^{-2}$},
mark size = 1,
height = 3.6cm,
width = 3.6cm]
\addplot[opacity=0.7, mark=o, mark size=1, blue] coordinates {(600,44.9055)(602,44.7205)(604,44.5235)(606,44.3132)(608,44.0881)(610,43.8467)(612,43.5872)(614,43.3073)(616,43.0046)(618,42.6762)(620,42.3186)(622,41.9279)(624,41.4991)(626,41.0264)(628,40.5027)(630,39.9192)(632,39.2652)(634,38.5269)(636,37.6869)(638,36.7228)(640,35.6047)(642,34.2924)(644,32.7308)(646,30.8412)(648,28.508)(650,25.5543)(652,21.6945)(654,16.4363)(656,8.85109)(658,-3.04463)(660,-24.3841)}; %q_1

\addplot[opacity=0.7, mark=x, mark size=1, red] coordinates {(600,52.8269)(602,52.6086)(604,52.3761)(606,52.1279)(608,51.8624)(610,51.5777)(612,51.2716)(614,50.9417)(616,50.5849)(618,50.1979)(620,49.7767)(622,49.3164)(624,48.8113)(626,48.2547)(628,47.6381)(630,46.9512)(632,46.1814)(634,45.3124)(636,44.324)(638,43.1895)(640,41.874)(642,40.3304)(644,38.4934)(646,36.2708)(648,33.5268)(650,30.0531)(652,25.5142)(654,19.3311)(656,10.412)(658,-3.57501)(660,-28.6653)}; %q_2
\end{axis}
\end{tikzpicture}
% Le = 10^{-1},
\begin{tikzpicture}
\pgfplotsset{every tick label/.append style={font=\tiny}, every axis/.append style={font=\tiny}}
\begin{axis}[
xlabel=$\R$,
ylabel = $q_i$,
title={$\Le=10^{-1}$},
scaled ticks=false,
axis lines=center,
y label style={at={(axis description cs:-0.1,.9)},anchor=south},
mark size = 1,
height = 3.6cm,
width = 3.6cm]
\addplot[opacity=0.7, mark=o, mark size=1, blue] coordinates {(600,2.52833)(603,2.46279)(606,2.39432)(609,2.32273)(612,2.24779)(615,2.16927)(618,2.08689)(621,2.00038)(624,1.9094)(627,1.81361)(630,1.71261)(633,1.60597)(636,1.49319)(639,1.37374)(642,1.247)(645,1.11229)(648,0.968821)(651,0.815718)(654,0.651975)(657,0.476443)(660,0.2878)(663,0.0845194)(666,-0.135171)(669,-0.373342)(672,-0.632428)(675,-0.915313)(678,-1.22543)(681,-1.56692)(684,-1.94479)(687,-2.36519)(690,-2.83572)(693,-3.36592)(696,-3.9679)(699,-4.65732)(702,-5.4547)(705,-6.38758)(708,-7.49369)};

\addplot[opacity=0.7, mark=x, mark size=1, red] coordinates {(600,2.97975)(603,2.90258)(606,2.82196)(609,2.73767)(612,2.64944)(615,2.557)(618,2.46003)(621,2.35819)(624,2.25111)(627,2.13837)(630,2.01951)(633,1.894)(636,1.76129)(639,1.62074)(642,1.47161)(645,1.3131)(648,1.14431)(651,0.964185)(654,0.771552)(657,0.565058)(660,0.343149)(663,0.104032)(666,-0.154379)(669,-0.434517)(672,-0.739247)(675,-1.07196)(678,-1.43669)(681,-1.8383)(684,-2.28268)(687,-2.77706)(690,-3.33039)(693,-3.95386)(696,-4.66174)(699,-5.47241)(702,-6.41002)(705,-7.50692)(708,-8.8075)};
\end{axis}
\end{tikzpicture}
% Le = 0.5,
\begin{tikzpicture}
\pgfplotsset{every tick label/.append style={font=\tiny}, every axis/.append style={font=\tiny}}
\begin{axis}[
xlabel=$\R$,
ylabel = $q_i$,
title={$\Le=0.5$},
scaled ticks=false,
axis lines=center,
y label style={at={(axis description cs:-0.1,1)},anchor=south},
mark size = 1,
height = 3.6cm,
width = 3.6cm]
\addplot[opacity=0.7, mark=o, mark size=1, blue] coordinates {(600,0.521524)(610,0.509056)(620,0.49627)(630,0.483154)(640,0.469693)(650,0.455876)(660,0.441687)(670,0.427111)(680,0.412133)(690,0.396734)(700,0.380898)(710,0.364605)(720,0.347835)(730,0.330567)(740,0.312778)(750,0.294445)(760,0.275542)(770,0.256042)(780,0.235916)(790,0.215133)(800,0.193662)(810,0.171467)(820,0.14851)(830,0.124753)(840,0.100151)(850,0.0746601)(860,0.0482303)(870,0.0208088)(880,-0.00766135)(890,-0.0372414)(900,-0.0679976)(910,-0.100002)(920,-0.133331)(930,-0.16807)(940,-0.204309)(950,-0.242149)(960,-0.281698)(970,-0.323074)(980,-0.366407)(990,-0.411839)(1000,-0.459527)(1010,-0.509642)(1020,-0.562376)(1030,-0.617939)(1040,-0.676564)(1050,-0.738513)(1060,-0.804077)(1070,-0.873581)(1080,-0.947391)(1090,-1.02592)(1100,-1.10964)};

\addplot[opacity=0.7, mark=x, mark size=1, red] coordinates {(600,0.618288)(610,0.603586)(620,0.58851)(630,0.573046)(640,0.557178)(650,0.54089)(660,0.524165)(670,0.506985)(680,0.489331)(690,0.471185)(700,0.452523)(710,0.433324)(720,0.413565)(730,0.393221)(740,0.372264)(750,0.350667)(760,0.3284)(770,0.305432)(780,0.281727)(790,0.257251)(800,0.231965)(810,0.205828)(820,0.178796)(830,0.150822)(840,0.121856)(850,0.0918446)(860,0.0607293)(870,0.0284482)(880,-0.00506559)(890,-0.0398842)(900,-0.0760856)(910,-0.113754)(920,-0.15298)(930,-0.193863)(940,-0.236511)(950,-0.28104)(960,-0.327578)(970,-0.376264)(980,-0.42725)(990,-0.480705)(1000,-0.53681)(1010,-0.595771)(1020,-0.657809)(1030,-0.723172)(1040,-0.792137)(1050,-0.865008)(1060,-0.942129)(1070,-1.02388)(1080,-1.1107)(1090,-1.20306)(1100,-1.30153)};
\end{axis}
\end{tikzpicture}
\begin{tikzpicture}
\pgfplotsset{every tick label/.append style={font=\tiny}, every axis/.append style={font=\tiny}}
\begin{axis}[
xlabel=$\R$,
ylabel = $q_i$,
title={$\Le=5$},
scaled ticks=false,
axis lines=center,
y label style={at={(axis description cs:-0.1,.9)},anchor=south},
mark size = 1,
height = 3.6cm,
width = 3.6cm]
\addplot[opacity=0.7, mark=o, mark size=1, blue] coordinates {(200,0.29079)(220,0.303498)(240,0.31522)(260,0.326065)(280,0.33613)(300,0.345495)(320,0.354231)(340,0.3624)(360,0.370055)(380,0.377243)(400,0.384007)(420,0.390382)(440,0.396402)(460,0.402095)(480,0.407488)(500,0.412603)(520,0.417462)(540,0.422084)(560,0.426485)(580,0.430682)(600,0.434687)(620,0.438515)(640,0.442176)(660,0.445682)(680,0.449042)(700,0.452265)(720,0.455359)(740,0.458333)(760,0.461192)(780,0.463944)(800,0.466595)(820,0.46915)(840,0.471614)(860,0.473992)(880,0.476288)(900,0.478507)(920,0.480653)(940,0.482729)(960,0.484739)(980,0.486685)(1000,0.488571)};

\addplot[opacity=0.7, mark=x, mark size=1, red] coordinates {(200,0.345274)(220,0.360275)(240,0.374117)(260,0.38693)(280,0.398825)(300,0.409898)(320,0.420233)(340,0.429901)(360,0.438965)(380,0.447482)(400,0.455499)(420,0.46306)(440,0.470203)(460,0.476963)(480,0.483369)(500,0.48945)(520,0.49523)(540,0.50073)(560,0.505972)(580,0.510973)(600,0.51575)(620,0.520318)(640,0.524691)(660,0.52888)(680,0.532899)(700,0.536756)(720,0.540463)(740,0.544027)(760,0.547458)(780,0.550762)(800,0.553947)(820,0.55702)(840,0.559985)(860,0.56285)(880,0.565619)(900,0.568298)(920,0.57089)(940,0.5734)(960,0.575832)(980,0.57819)(1000,0.580477)};
\end{axis}
\end{tikzpicture}
\caption{The plots of $q_1$(blue with circled mark) and $q_2$ (red with x mark) against $\R$ for $\Le=10^{-2},\ 10^{-1},\ 0.5,\ 5$ with $r=2/\pi$ for $l_c=1$ and $r=2\sqrt{3}/\pi$ for $l_c=2$ and $\Pre=7.5$. The type of transition changes from a Type-I transition to a Type-II transition at the $\R$ value corresponding to a change of sign of $q_{l_c}$.}
\label{fig:q1 q2 vs R}
\end{figure}
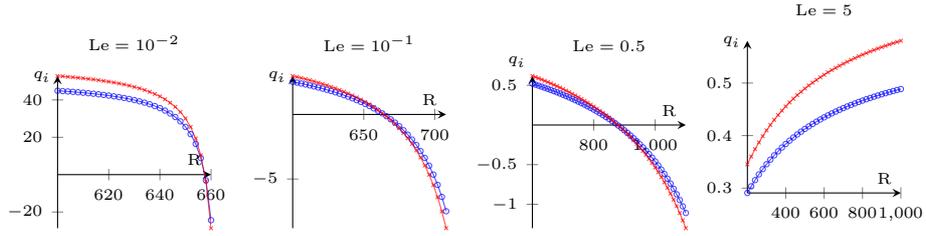

For $\Le=10^{-2},$ $\Le=10^{-1}$, $\Le=0.5$ and in both $l_c=1$ and $l_c=2$ cases, we find that the transition is Type-I ($q_{l_c} > 0$) if $\R < \R_{\ast}$ and Type-II ($q_{l_c} < 0$) if $\R_{\ast} < \R < \R_0$ as shown in Table~\ref{lc=1 results}, Table~\ref{lc=2 results} and Figure~\ref{fig:q1 q2 vs R} with $R_{\ast} \approx R_1 = \frac{\sigma_c}{1-\Le^2}$. The numerical values of $\R_1$ when $\sigma_c = \frac{27}{4} \pi^4$ is $\R_1 = 657.577$ when $\Le=10^{-2}$, $\R_1 = 664.153$ when $\Le=10^{-1}$ and $\R_1 = 876.682$ when $\Le = 0.5$. For $\Le=5$, $q_{l_c}$ does not change sign and the transition is Type-I ($q_{l_c} > 0$) for all $\R$, see Figure~\ref{fig:q1 q2 vs R}.
% N[((27/4)*Pi^4)/(1 - Le^2) /. Le -> {10^(-2), 10^(-1), 0.5}]
% Out = {657.577, 664.153, 876.682}

\begin{table}[ht]
\centering
\caption{Types of transitions with respect to $\R$ for $l_c=1$, $r=2/\pi$, $\Pre=7.5 $ for different $\Le$ values. \label{lc=1 results}}
\begin{tabular}{|l|l|l|}
\hline
$\Le $  & Type-I transition   & Type-II transition \\ \hline
$10^{-2}$ & if $\R < 657.577$ & if $657.577< \R < 665.038$ \\ \hline
$10^{-1}$ & if $\R < 664.182$ & if $664.182 < \R < 740.309$ \\ \hline
$0.5$ & if $\R < 877.346$ & if $877.346 < \R < 1402.69$ \\ \hline
\end{tabular}
\end{table}

\begin{table}[ht]
\centering
\caption{Types of transitions with respect to $\R$ for $l_c=2$, $r=2\sqrt{3}/\pi$, $\Pre=7.5 $ for different $\Le$ values.
\label{lc=2 results}}
\begin{tabular}{|l|l|l|}
\hline
$\Le $  & Type-I transition   & Type-II transition \\ \hline
$10^{-2}$ & if $\R < 657.578$ & if $657.578< \R < 665.038$ \\ \hline
$10^{-1}$ & if $\R < 664.236$ & if $664.236 < \R < 740.309$ \\ \hline
$0.5$ & if $\R < 878.513$ & if $878.513 < \R < 1402.69$ \\ \hline
\end{tabular}
\end{table}

% \bibliography{thermo.bib}{}
\bibliographystyle{siam}

\end{document}